\documentclass[small]{article}
\usepackage{amsmath,amstext,amsthm,amsfonts}
\usepackage{makeidx}
\usepackage{amssymb,latexsym}
\usepackage{txfonts, pxfonts}
\usepackage[english]{babel}
\usepackage[pdftex]{graphicx}
\numberwithin{equation}{section}
\newtheorem{theorem}{\textbf{Theorem}}[section]

\newtheorem{lemma}[theorem]{\textbf{Lemma}}

\theoremstyle{definition}
\newtheorem{definition}[theorem]{\textbf{Definition}}
\theoremstyle{notation}
\newtheorem*{notation}{\textbf{Notation}}
\theoremstyle{remark}
\newtheorem{remark}[theorem]{\it{Remark}}

\def\A{\mathcal{A}}
\def\e{\epsilon}

\def\R{\mathbb{R}}
\def\Rn{{\mathbb{R}}^n_+}
\def\d{\partial}

\def\dv{dv}
\def\ds{d\sigma}

\def\b{\beta}
\def\l{\lambda}

\def\tg{\tilde{g}_{\nu}}

\def\critbordo{\frac{2(n-1)}{n-2}}

\def\conjbordo{\frac{2(n-1)}{n}}\def\ba{\begin{align}}

\def\u{\hat{u}_{\nu}}

\def\conjbordo{\frac{2(n-1)}{n}}\def\ba{\begin{align}}
\def\ea{\end{align}}
\def\bp{\begin{proof}}
\def\ep{\end{proof}}

\title{The asymptotic behavior of Palais-Smale sequences on manifolds with boundary}
\author{ \textsc{S\'ergio Almaraz}~\footnote{Supported by CAPES and FAPERJ (Brazil).}}
\date{}

\begin{document}
\maketitle

\begin{abstract}
We describe the asymptotic behavior of Palais-Smale sequences associated to certain Yamabe-type equations on manifolds with boundary. 
We prove that each of those sequences converges to a solution of the limit equation plus a finite number of "bubbles" which are  obtained by rescaling fundamental solutions of the corresponding  Euclidean equations.
\end{abstract}

\section{Introduction}

Let $(M^n,g)$ be a compact Riemannian manifold with boundary $\d M$ and dimension $n\geq 3$. For $u\in H^1(M)$, we consider the following family of equations, indexed by $\nu \in\mathbb{N}$:
\begin{equation}\label{eq:Enu}
\begin{cases}
\Delta_g u = 0\,,&\text{in}\:M\,,\\
\frac{\d}{\d \eta_g}u-h_{\nu}u+u^{\frac{n}{n-2}}=0\,,&\text{on}\:\d M\,,
\end{cases}
\end{equation}
and their associated functionals
\begin{equation}\label{def:I}
I_g^{\nu}(u)=\frac{1}{2}\int_{M}|du|_g^2\dv_g+\frac{1}{2}\int_{\d M}h_{\nu}u^2\ds_g
-\frac{n-2}{2(n-1)}\int_{\d M}|u|^{\critbordo}\ds_g\,.
\end{equation}
Here, $\{h_{\nu}\}_{\nu\in\mathbb{N}}$ is a sequence of  functions in $C^{\infty}(\d M)$, $\Delta_g$ is the Laplace-Beltrami operator, and $\eta_g$ is the inward unit normal vector to $\d M$. Moreover, $\dv_g$ and $\ds_g$ are the volume forms of $M$ and $\d M$ respectively and $H^1(M)$ is the Sobolev space $H^1(M)=\{u\in L^2(M)\,;\:du\in L^2(M)\}$.
\begin{definition}
We say that $\{u_{\nu}\}_{\nu\in\mathbb{N}}\subset H^1(M)$ is a {\it{Palais-Smale}} sequence for $\{I_g^{\nu}\}$ if\\
(i) $\{I_g^{\nu}(u_{\nu})\}_{\nu\in\mathbb{N}}$ is bounded, and\\
(ii) $dI_g^{\nu}(u_{\nu})\to 0$ strongly in $H^1(M)'$ as $\nu\to\infty$.
\end{definition}

In this paper we establish a result describing the asymptotic behavior of those  Palais-Smale sequences.
This work is inspired by Struwe's theorem  in \cite{struwe} for equations $\Delta u+\l u+|u|^{\frac{4}{n-2}}u=0$ on Euclidean domains. We refer the reader to \cite[Chapter 3]{druet-hebey-robert} for a version of Struwe's theorem on closed Riemannian manifolds, and to \cite{cao-noussair-yan, chabrowski-girao, pierotti-terracini} for similar equations with boundary conditions.

Roughly speaking, as $\nu\to\infty$ and $h_{\nu}\to h_{\infty}$ we prove that each  Palais-Smale sequence $\{u_{\nu}\geq 0\}_{\nu\in\mathbb{N}}$  is $H^1(M)$-asymptotic to a nonnegative solution of the limit equations 
\begin{equation}\label{eq:Einfty}
\begin{cases}
\Delta_g u = 0\,,&\text{in}\:M\,,\\
\frac{\d}{\d \eta_g}u-h_{\infty}u+u^{\frac{n}{n-2}}=0\,,&\text{on}\:\d M\,,
\end{cases}
\end{equation}
plus a finite number of "bubbles" obtained by rescaling fundamental positive solutions of the Euclidean equations
\begin{equation}
\label{eq:U}
\begin{cases}
\Delta u = 0\,,&\text{in}\:\mathbb{R}_+^n\,,\\
\frac{\partial}{\partial y_n}u+ u^{\frac{n}{n-2}}=0\,,&\text{on}\:\partial\mathbb{R}_+^n\,,
\end{cases}
\end{equation}
where $\Rn=\{(y_1,...,y_n)\in \R^n\,;\:y_n\geq 0\}$.

Palais-Smale sequences frequently appear in the blow-up analysis of geometric problems. In the particular case when $h_{\infty}$ is $\frac{n-2}{2}$ times the boundary mean curvature, the equations (\ref{eq:Einfty}) are satisfied by a positive smooth function $u$ representing a conformal scalar-flat Riemannian metric $u^{\frac{4}{n-2}}g$ with positive constant boundary mean curvature. The existence of those metrics is the Yamabe-type problem for manifolds with boundary introduced by Escobar in \cite{escobar3}.

An application of our result is the blow-up analysis performed by the author in \cite{almaraz4} for the proof of a convergence theorem for a Yamabe-type flow  introduced by Brendle in \cite{brendle-boundary}.

We now begin to state our theorem more precisely.
 
\vspace{0.2cm}\noindent
{\bf{Convention.}} We assume that there is some $h_{\infty}\in C^{\infty}(\d M)$ and $C>0$ such that $h_{\nu}\to h_{\infty}$ in $L^2(\d M)$ as $\nu\to\infty$ and $|h_{\nu}(x)|\leq C$ for all $x\in\d M$, $\nu\in\mathbb{N}$. This obviously implies that $h_{\nu}\to h_{\infty}$ in $L^p(\d M)$ as $\nu\to\infty$, for any $p\geq 1$.

\begin{notation}
If $(M,g)$ is a Riemannian manifold with boundary $\d M$, 
we will denote by $D_r(x)$ the metric ball in $\d M$ with center at $x\in \d M$ and radius $r$.

If $z_0\in \Rn$, we set $B_r^+(z_0)=\{z\in\Rn\,;\:|z-z_0|< r\}$. We define 
$$
\d^+B^+_{r}(z_0)=\d B^+_{r}(z_0)\cap\Rn\,,
\:\:\:\:\text{and}\:\:\:\:\:
\d 'B^+_{r}(z_0)=B^+_{r}(z_0)\cap \d\Rn\,.
$$
Thus, $\d 'B^+_{r}(z_0)=\emptyset$ if $z_0=(z_0^1,...,z_0^n)$ satisfies $z_0^n> r$.
\end{notation}

We define the Sobolev space $D^1(\Rn)$ as the completion of $C_0^{\infty}(\Rn)$ with respect to the norm
$$
\|u\|_{D^1(\Rn)}=\sqrt{\int_{\Rn}|du(y)|^2dy}\,. 
$$

It follows from a Liouville-type theorem established by Li and Zhu in \cite{li-zhu} (see also \cite{escobar1} and \cite{chipot-shafrir-fila}) that any nonnegative solution in $D^1(\Rn)$ to the equations (\ref{eq:U}) is of the form 
\begin{equation}\label{def:U}
U_{\e,a}(y)=\left(\frac{\e}{(y_n+\frac{\e}{n-2})^2+|\bar{y}-a|^2}\right)^{\frac{n-2}{2}}\,,
\:\:\:\:\:\:a\in\R^{n-1},\,\e>0\,,
\end{equation}
or is identically zero (see Remark \ref{rm:reg}). By Escobar (\cite{escobar0}) or Beckner (\cite{beckner}) we have the sharp Euclidean Sobolev inequality
\begin{equation}\label{des:sobolev}
\left(\int_{\d\Rn}|u(y)|^{\critbordo}dy\right)^{\frac{n-2}{n-1}}\leq 
K_n^2\int_{\Rn}|du(y)|^2dy\,,
\end{equation}
for $u\in D^1(\Rn)$,
which has the family of functions (\ref{def:U}) as extremal functions. Here, 
$$
K_n=\left(\frac{n-2}{2}\right)^{-\frac{1}{2}}\sigma_{n-1}^{-\frac{1}{2(n-1)}}\,,
$$
where $\sigma_{n-1}$ is the area of the unit $ (n-1)$-sphere in $\R^n$. Up to a multiplicative  constant, the functions defined by (\ref{def:U}) are the only nontrivial extremal ones for  the inequality (\ref{des:sobolev}).
\begin{definition}
Fix $x_0\in\d M$ and geodesic normal coordinates for $\d M$ centered at $x_0$. 
Let $(x_1,...,x_{n-1})$ be the coordinates of  $x\in\d M$ and $\eta_g(x)$ be the inward unit vector normal to $\d M$ at $x$. 
For small $x_n\geq 0$, the point $\exp_{x}(x_n\eta_g(x))\in M$ is said to have {\it{Fermi coordinates}} $(x_1,...,x_n)$ (centered at $x_0$). 
\end{definition}
For small $\rho>0$ the Fermi coordinates centered at $x_0\in \d M$ define a smooth map $\psi_{x_0}:B^+_{\rho}(0)\subset \Rn\to M$.

We define the functional $I_g^{\infty}$  by the same expression as $I_g^{\nu}$ with $h_{\nu}=h_{\infty}$ for all $\nu$, and state our main theorem as follows:
\begin{theorem}\label{struwe_thm}
Let $(M^n,g)$ be a compact  Riemannian manifold with boundary $\d M$ and dimension $n\geq 3$. Suppose $\{u_{\nu}\geq 0\}_{\nu\in\mathbb{N}}$ is a Palais-Smale sequence for $\{I_g^{\nu}\}$. Then there exist $m\in\{0,1,2,...\}$, a nonnegative solution  $u^0\in H^1(M)$ of (\ref{eq:Einfty}), and $m$ nontrivial nonegative solutions $U^j=U_{\e_j,a_j}\in D^1(\Rn)$ of (\ref{eq:U}), sequences  $\{R_{\nu}^j>0\}_{\nu\in\mathbb{N}}$, and sequences $\{x_{\nu}^j\}_{\nu\in\mathbb{N}}\subset \d M$, $1\leq j\leq m$, the whole satisfying the following conditions for $1\leq j\leq m$, possibly after taking subsequences: 
\\
(i) $R_{\nu}^j\to\infty$ as $\nu\to\infty$.
\\
(ii) $x_{\nu}^j$ converges as $\nu\to\infty$.
\\
(iii) $\big\| u_{\nu}-u^0-\sum_{j=1}^{m}\eta_{\nu}^ju_{\nu}^j\big\|_{H^1(M)}\to 0$ as $\nu\to\infty$, where 
$$
u_{\nu}^j(x)=(R_{\nu}^j)^{\frac{n-2}{2}}U^j(R_{\nu}^j\psi^{-1}_{x_{\nu}^j}(x))
\:\:\:\:\:\text{for}\:x\in \psi_{x_{\nu}^j}(B^+_{2r_0}(0))\,.
$$
Here, $r_0>0$ is small, the
$$\psi_{x_{\nu}^j}:B^+_{2r_0}(0)\subset \Rn\to M$$
are Fermi coordinates centered at $x_{\nu}^j\in \d M$, and the $\eta_{\nu}^j$ are smooth cutoff functions such that $\eta_{\nu}^j\equiv 1$ in $\psi_{x_{\nu}^j}(B^+_{r_0}(0))$ and $\eta_{\nu}^j\equiv 0$ in $M\backslash\psi_{x_{\nu}^j}(B^+_{2r_0}(0))$.

Moreover, 
$$
I_g^{\nu}(u_{\nu})-I_g^{\infty}(u^0)-\frac{m}{2(n-1)}K_n^{-2(n-1)}\to 0
\:\:\:\:\text{as}\:\nu\to\infty\,,
$$
and we can assume that for all $i\neq j$
\begin{equation}\label{eq:brezis-coron}
\frac{R_{\nu}^i}{R_{\nu}^j}+\frac{R_{\nu}^j}{R_{\nu}^i}+R_{\nu}^iR_{\nu}^jd_g(x_{\nu}^i,x_{\nu}^j)^2\to \infty
\:\:\:\:\text{as}\:\nu\to\infty\,.
\end{equation}
\end{theorem}
\begin{remark}
Relations of the type (\ref{eq:brezis-coron}) were previously obtained in \cite{bahri-coron, brezis-coron}.
\end{remark}

\section{Proof of the main theorem}

The rest of this paper is devoted to the proof of Theorem \ref{struwe_thm} which will be carried out in several lemmas. Our presentation will follow the same steps as Chapter 3 of \cite{druet-hebey-robert}, with the necessary modifications.
\begin{lemma}\label{step1}
Let $\{u_{\nu}\}$ be a Palais-Smale sequence for $\{I_g^{\nu}\}$. Then there exists $C>0$ such that $\|u_{\nu}\|_{H^1(M)}\leq C$ for all $\nu$.
\end{lemma}
\bp
It suffices to prove that $\|du_{\nu}\|_{L^2(M)}$ and $\|u_{\nu}\|_{L^2(\d M)}$ are uniformly bounded. The proof follows the same arguments as \cite[p.27]{druet-hebey-robert}.
\ep

Define $I_g$ as the functional $I_g^{\nu}$ when $h_{\nu}\equiv 0$ for all $\nu$.
\begin{lemma}\label{step2}
Let $\{u_{\nu}\geq 0\}$ be a Palais-Smale sequence for $\{I_g^{\nu}\}$ such that $u_{\nu}\rightharpoonup u^0\geq 0$ in $H^1(M)$ and set $\hat{u}_{\nu}=u_{\nu}-u^0$. Then $\{\hat{u}_{\nu}\}$ is a Palais-Smale sequence for $\{I_g\}$ and satisfies
\begin{equation}\label{step2:1}
I_g(\hat{u}_{\nu})-I_g^{\nu}(u_{\nu})+I_g^{\infty}(u^0)\to 0
\:\:\:\:\text{as}\:\nu\to\infty\,.
\end{equation}
Moreover, $u^0$ is a (weak) solution of (\ref{eq:Einfty}).
\end{lemma}
\bp
First observe that $u_{\nu}\rightharpoonup u^0$ in $H^1(M)$ implies that $u_{\nu}\to u^0$ in $L^{\frac{n}{n-2}}(\d M)$ and a.e. in $\d M$. Using the facts that $dI_g^{\nu}(u_{\nu})\phi\to 0$ for any $\phi\in C^{\infty}(\bar{M})$ and $h_{\nu}\to h_{\infty}$ in $L^p(\d M)$ for any $p\geq 1$, it is not difficult to see that the last assertion of Lemma \ref{step2} follows.

In order to prove (\ref{step2:1}), we first observe that 
$$
I_g^{\nu}(u_{\nu})=I_g(\hat{u}_{\nu})+I_g^{\infty}(u^0)-\frac{(n-2)}{2(n-1)}\int_{\d M}\Phi_{\nu}\ds_g+o(1)\,,
$$
where $\Phi_{\nu}=|\hat{u}_{\nu}+u^0|^{\critbordo}-|\hat{u}_{\nu}|^{\critbordo}-|u^0|^{\critbordo}$, and $o(1)\to 0$ as $\nu\to \infty$. 
Then (\ref{step2:1}) follows from the fact that there exists $C>0$ such that
$$
\int_{\d M}\Phi_{\nu}\ds_g\leq C\int_{\d M}|\hat{u}_{\nu}|^{\frac{n}{n-2}}|u^0|\ds_g
+C\int_{\d M}|u^0|^{\frac{n}{n-2}}|\hat{u}_{\nu}|\ds_g\,,
\:\:\:\:\text{for all}\:\nu\,,
$$
and, by basic integration theory, the right side of this last inequality goes to $0$ as $\nu\to\infty$.

Now we prove that $\{\hat{u}_{\nu}\}$ is a Palais-Smale sequence for $I_g$. Let $\phi\in C^{\infty}(M)$. Observe that
$$
\left|\int_{\d M}h_{\nu}u_{\nu}\phi\ds_g-\int_{\d M}h_{\infty}u_{\nu}\phi\ds_g\right|
\leq 
\|u_{\nu}\|_{L^2(\d M)}\|h_{\nu}-h_{\infty}\|_{L^{2(n-1)}(\d M)}\|\phi\|_{L^{\critbordo}(\d M)}
$$
by H\"{o}lder's inequality. Then, by the Sobolev embedding theorem,
$$
\int_{\d M}h_{\nu}u_{\nu}\phi\ds_g=\int_{\d M}h_{\infty}u^0\phi\ds_g+o(\|\phi\|_{H^1(M)})
$$ 
from which follows that
\begin{equation}\label{step2:2}
dI_g^{\nu}(u_{\nu})\phi=dI_g(\hat{u}_{\nu})\phi
-\int_{\d M}\psi_{\nu}\phi\ds_g+o(\|\phi\|_{H^1(M)})\,,
\end{equation}
where $\psi_{\nu}=|\hat{u}_{\nu}+u^0|^{\frac{2}{n-2}}(\hat{u}_{\nu}+u^0)
-|\hat{u}_{\nu}|^{\frac{2}{n-2}}\hat{u}_{\nu}-|u^0|^{\frac{2}{n-2}}u^0$.

Next we observe that there exists $C>0$ such that
$$
\int_{\d M}|\psi_{\nu}\phi|\ds_g
\leq C\int_{\d M}|\hat{u}_{\nu}|^{\frac{2}{n-2}}|u^0|\,|\phi|\ds_g
+C\int_{\d M}|u^0|^{\frac{2}{n-2}}|\hat{u}_{\nu}|\,|\phi|\ds_g\,,
$$
for all $\nu$, and use H\"{o}lder's inequality and basic integration theory to obtain
\ba
\int_{\d M}|\psi_{\nu}\phi|\ds_g
&\leq
\big(\big{\|}|\hat{u}_{\nu}|^{\frac{2}{n-2}}u^0\big{\|}_{L^{\conjbordo}(\d M)}
+\big{\|}|u^0|^{\frac{2}{n-2}}\hat{u}_{\nu}\big{\|}_{L^{\conjbordo}(\d M)}\big)\,
\|\phi\|_{L^{\critbordo}(\d M)}\notag
\\
&=o\big(\|\phi\|_{L^{\critbordo}(\d M)}\big)\,.\notag
\end{align}
Then we can use this and the Sobolev embedding theorem in (\ref{step2:2}) to conclude that
$$
dI_g^{\nu}(u_{\nu})\phi=dI_g(\hat{u}_{\nu})\phi
+o(\|\phi\|_{H^1(M)})\,,
$$
finishing the proof.
\ep
\begin{lemma}\label{step3}
Let $\{\hat{u}_{\nu}\}_{\nu\in\mathbb{N}}$ be a Palais-Smale sequence for $I_g$ such that $\hat{u}_{\nu}\rightharpoonup 0$ in $H^1(M)$ and $I_g(\hat{u}_{\nu})\to \b$ as $\nu\to\infty$ for some $\displaystyle{\b<\frac{K_n^{-2(n-1)}}{2(n-1)}}$. Then $\hat{u}_{\nu}\to 0$ in $H^1(M)$ as $\nu\to\infty$.
\end{lemma}
\bp
Since 
$$
\int_M|d\hat{u}_{\nu}|^2\dv_g-\int_{\d M}|\hat{u}_{\nu}|^{\critbordo}\ds_g
=dI_g(\hat{u}_{\nu})\cdot \hat{u}_{\nu}
=o(\|\hat{u}_{\nu}\|_{H^1(M)})
$$
and $\{\|\hat{u}_{\nu}\|_{H^1(M)}\}$ is uniformly bounded due to Lemma \ref{step1}, we can see that
\ba\label{step1:1}
\b+o(1)
=I_g(\hat{u}_{\nu})
&=\frac{1}{2(n-1)}\int_{\d M}|\hat{u}_{\nu}|^{\critbordo}\ds_g+o(1)
\\
&=\frac{1}{2(n-1)}\int_{M}|d\hat{u}_{\nu}|_g^2\dv_g+o(1)\notag
\end{align}
which already implies $\b\geq 0$. At the same time, as proved by Li and Zhu in \cite{li-zhu2}, there exists $B=B(M,g
)>0$ such that
$$
\left(\int_{\d M}|\hat{u}_{\nu}|^{\critbordo}\ds_g\right)^{\frac{n-2}{n-1}}
\leq
K_n^2\int_{M}|d\hat{u}_{\nu}|_g^2\dv_g+B\int_{\d M}|\hat{u}_{\nu}|^2\ds_g\,.
$$
Since $H^1(M)$ is compactly embedded in $L^2(\d M)$, we have $\|\hat{u}_{\nu}\|_{L^2(\d M)}\to 0$. Then we obtain
$$
(2(n-1)\b+o(1))^{\frac{n-2}{n-1}}
\leq
2(n-1)K_n^2\b+o(1)
$$
from which we conclude that either $$\frac{K_n^{-2(n-1)}}{2(n-1)}\leq \b + o(1)$$ or $\b=0$. Hence, our hypotheses imply $\b=0$. Using (\ref{step1:1}) finishes the proof.  
\ep

Define the functional
$$
E(u)=\frac{1}{2}\int_{\Rn}|du(y)|^2dy-\frac{n-2}{2(n-1)}\int_{\d\Rn}|u(y)|^{\critbordo}dy
$$
for $u\in D^1(\Rn)$ and observe that $E(U_{\e,a})={\displaystyle{\frac{K_n^{-2(n-1)}}{2(n-1)}}}$ for any $a\in\R^{n-1}$, $\e>0$. 
\begin{lemma}\label{step4}
Let $\{\hat{u}_{\nu}\}_{\nu\in\mathbb{N}}$ be a Palais-Smale sequence for $I_g$. Suppose $\hat{u}_{\nu}\rightharpoonup 0$ in $H^1(M)$, but not strongly. Then there exist a sequence $\{R_{\nu}>0\}_{\nu\in\mathbb{N}}$ with $R_{\nu}\to\infty$, a convergent sequence $\{x_{\nu}\}_{\nu\in\mathbb{N}}\subset \d M$, and a nontrivial solution $u\in D^1(\Rn)$ of 
\begin{equation}\label{eq:u}
\begin{cases}
\Delta u=0\,,&\text{in}\:\Rn\,,
\\
\frac{\d}{\d y_n}u-|u|^{\frac{2}{n-2}}u=0\,,&\text{on}\:\d\Rn\,,
\end{cases}
\end{equation}
the whole such that, up to a subsequence, the following holds: If 
$$
\hat{v}_{\nu}(x)=\hat{u}_{\nu}(x)-\eta_{\nu}(x)R_{\nu}^{\frac{n-2}{2}}u(R_{\nu}\psi_{x_{\nu}}^{-1}(x))\,,
$$ 
then $\{\hat{v}_{\nu}\}_{\nu\in\mathbb{N}}$ is a Palais-Smale sequence for $I_g$ satisfying $\hat{v}_{\nu}\rightharpoonup 0$ in $H^1(M)$ and 
$$
\lim_{\nu\to\infty}\big(I_g(\hat{u}_{\nu})-I_g(\hat{v}_{\nu})\big)=E(u)\,.
$$
Here, the $\psi_{x_{\nu}}:B^+_{2r_0}(0)\subset\Rn\to M$
are Fermi coordinates centered at $x_{\nu}$
and the $\eta_{\nu}(x)$ are smooth cutoff functions such that $\eta_{\nu}\equiv 1$ in $\psi_{x_{\nu}}(B^+_{r_0}(0))$ and $\eta_{\nu}\equiv 0$ in $M\backslash \psi_{x_{\nu}}(B^+_{2r_0}(0))$.
\end{lemma}
\bp
By the density of $C^{\infty}(M)$ in $H^1(M)$ we can assume that $\hat{u}_{\nu}\in C^{\infty}(M)$. We can also assume that $I_g(\hat{u}_{\nu})\to\b$ as $\nu\to\infty$ and, since $dI_g(\hat{u}_{\nu})\to 0$ in $H^1(M)'$, we obtain
$$
\lim_{\nu\to\infty}\int_{\d M} |\hat{u}_{\nu}|^{\critbordo}\ds_g=2(n-1)\b\geq K_n^{-2(n-1)}
$$
as in the proof of Lemma \ref{step3}.
Hence, given $t_0>0$ small we can choose $x_0\in\d M$ and $\l_0>0$ such that
$$
\int_{D_{t_0}(x_0)}|\hat{u}_{\nu}|^{\critbordo}\ds_g\geq \l_0
$$
up to a subsequence. Now we set 
$$
\mu_{\nu}(t)=\max_{x\in\d M}\int_{D_t(x)}|\hat{u}_{\nu}|^{\critbordo}\ds_g
$$ 
for $t>0$, and, for any $\l\in (0,\l_0)$, choose sequences $\{t_{\nu}\}\subset (0,t_0)$ and $\{x_{\nu}\}\subset\d M$ such that
\begin{equation}\label{def:tnu}
\l=\mu_{\nu}(t_{\nu})=\int_{D_{t_{\nu}}(x_{\nu})}|\hat{u}_{\nu}|^{\critbordo}\ds_g\,.
\end{equation} 
We can also assume that $x_{\nu}$ converges. Now we choose $r_0>0$ small such that for any $x_0\in \d M$ the Fermi coordinates $\psi_{x_0}(z)$ centered at $x_0$ are defined for all $z\in B^+_{2r_0}(0)\subset\Rn$ and satisfy
$$
\frac{1}{2}|z-z'|\leq d_g(\psi_{x_0}(z),\psi_{x_0}(z'))\leq 2|z-z'|\,,
\:\:\:\: 
\text{for any}\:z,z'\in B^+_{r_0}(0)\,.
$$

For each $\nu$ we consider Fermi coordinates 
$$
\psi_{\nu}=\psi_{x_{\nu}}:B^+_{2r_0}(0)\to M\,.
$$
For any $R_{\nu}\geq 1$ and $y\in B^+_{R_{\nu}r_0}(0)$ we set 
$$
\tilde{u}_{\nu}(y)=R_{\nu}^{-\frac{n-2}{2}}\hat{u}_{\nu}(\psi_{\nu}(R_{\nu}^{-1}y))
\:\:\:\:\text{and}\:\:\:\:
\tg(y)=(\psi_{\nu}^*g)(R_{\nu}^{-1}y)\,.
$$
Let us consider $z\in\Rn$ and $r>0$ such that $|z|+r<R_{\nu}r_0$.
Then we have 
\begin{equation}\notag
\int_{B^+_r(z)}|d\tilde{u}_{\nu}|_{\tg}^2\dv_{\tg}
=\int_{\psi_{\nu}(R_{\nu}^{-1}B^+_r(z))}|d\hat{u}_{\nu}|_g^2\dv_g\,,
\end{equation}
and, if  in addition $z\in\d\Rn$,
\begin{align}\label{step4:1'}
\int_{\d 'B^+_r(z)}|\tilde{u}_{\nu}|^{\critbordo}\ds_{\tg}
&=\int_{\psi_{\nu}(R_{\nu}^{-1}\d 'B^+_r(z))}|\hat{u}_{\nu}|^{\critbordo}\ds_g
\\
&\leq
\int_{D_{2R_{\nu}^{-1}r}(\psi_{\nu}(R_{\nu}^{-1}z))}|\hat{u}_{\nu}|^{\critbordo}\ds_g \,,\notag
\end{align}
where we have used the fact that
\begin{equation}\notag
\psi_{\nu}(R_{\nu}^{-1}\d 'B^+_r(z))
=\psi_{\nu}(\d 'B^+_{R_{\nu}^{-1}r}(R_{\nu}^{-1}z))
\subset D_{2R_{\nu}^{-1}r}(\psi_{\nu}(R_{\nu}^{-1}z))\,.
\end{equation}

Given $r\in(0,r_0)$ we fix $t_0\leq 2r$. Then, given a $\l\in(0,\l_0)$ to be fixed later, we set $R_{\nu}=2rt_{\nu}^{-1}\geq 2rt_0^{-1}\geq 1$. Then it follows from (\ref{def:tnu}) and (\ref{step4:1'}) that 
\begin{equation}\label{step4:3}
\int_{\d 'B^+_r(z)}|\tilde{u}_{\nu}|^{\critbordo}\ds_{\tg}\leq \l\,.
\end{equation}
Moreover, since $\psi_{\nu}(\d 'B^+_{2R_{\nu}^{-1}r}(0))=D_{t_{\nu}}(x_{\nu})$ we have
\begin{equation}\label{step4:4}
\int_{\d 'B^+_{2r}(0)}|\tilde{u}_{\nu}|^{\critbordo}\ds_{\tg}
=
\int_{D_{t_{\nu}}(x_{\nu})}|\hat{u}_{\nu}|^{\critbordo}\ds_{g}
=\l\,.
\end{equation}

Choosing $r_0$ smaller if necessary, we can suppose that
\begin{equation}\label{step4:5}
\frac{1}{2}\int_{\Rn}|du|^2dy
\leq
\int_{\Rn}|du|_{\tilde{g}_{x_0,R}}^2\dv_{\tilde{g}_{x_0,R}}
\leq 
2\int_{\Rn}|du|^2dy
\end{equation}
for any $R\geq 1$ and any $u\in D^1(\Rn)$ such that $\text{supp}(u)\subset B^+_{2r_0R}(0)$. Here, $\tilde{g}_{x_0,R}(y)=(\psi_{x_0}^*g)(R^{-1}y)$. We can also assume that 
\begin{equation}\label{step4:5'}
\frac{1}{2}\int_{\d\Rn}|u|dy
\leq
\int_{\d\Rn}|u|\ds_{\tilde{g}_{x_0,R}}
\leq 
2\int_{\d\Rn}|u|dy
\end{equation}
for all $u\in L^1(\d\Rn)$ such that $\text{supp}(u)\subset \d 'B^+_{2r_0R}(0)$.

Let $\tilde{\eta}$ be a smooth cutoff function on $\R^n$ such that $0\leq \tilde{\eta}\leq 1$, $\tilde{\eta}(z)= 1$ for $|z|\leq \frac{1}{4}$, and $\tilde{\eta}(z)= 0$ for $|z|\geq \frac{3}{4}$ . We set $\tilde{\eta}_{\nu}(y)=\tilde{\eta}(r_0^{-1}R_{\nu}^{-1}y)$.

It is easy to check that $\big{\{}\int_{\Rn}|d(\tilde{\eta}_{\nu}\tilde{u}_{\nu})|_{\tg}^2\dv_{\tg}\big{\}}$ is uniformly bounded. Then the inequality (\ref{step4:5}) implies that $\{\tilde{\eta}_{\nu}\tilde{u}_{\nu}\}$ is uniformly bounded in $D^1(\Rn)$ and we can assume that $\tilde{\eta}_{\nu}\tilde{u}_{\nu}\rightharpoonup u$ in $D^1(\Rn)$ for some $u$.

\bigskip\noindent
{\it{Claim 1.}} Let us set $\displaystyle r_1=r_0/24$. There exists $\l_1=\l_1(n)$ such that for any $0<r<r_1$ and $0<\l<\min\{\l_1,\l_0\}$ we have
$$
\tilde{\eta}_{\nu}\tilde{u}_{\nu}\to u\,,
\:\:\:\:\text{in}\:H^1(B^+_{2Rr}(0))\,,
\:\:\:\:\text{as}\:\nu\to\infty\,,
$$ 
for any $R\geq 1$ satisfying $R\leq R_{\nu}$ for all $\nu$ large.

\vspace{0.2cm}\noindent
{\it{Proof of Claim 1.}} We consider $r\in (0,r_1)$, $\l\in (0,\l_0)$  and choose $z_0\in\d\Rn$ such that $|z_0|<3(2R-1)r_1$. By Fatou's lemma,
\ba
\int_r^{2r}\liminf_{\nu\to\infty}
&\left\{
\int_{\d^+B^+_{\rho}(z_0)}\Big\{|d(\tilde{\eta}_{\nu}\tilde{u}_{\nu})|^2+|\tilde{\eta}_{\nu}\tilde{u}_{\nu}|^2\Big\}\ds_{\rho}
\right\}d\rho\notag
\\
&\leq
\liminf_{\nu\to\infty}
\int_{B^+_{2r}(z_0)}\Big\{|d(\tilde{\eta}_{\nu}\tilde{u}_{\nu})|^2+|\tilde{\eta}_{\nu}\tilde{u}_{\nu}|^2\Big\}dy
\leq C\,,\notag
\end{align}
where $\ds_{\rho}$ is the volume form on $\d^+B^+_{\rho}(z_0)$ induced by the Euclidean metric. 
Thus there exists $\rho\in [r,2r]$ such that, up to a subsequence,
$$
\int_{\d^+B^+_{\rho}(z_0)}\Big\{|d(\tilde{\eta}_{\nu}\tilde{u}_{\nu})|^2+|\tilde{\eta}_{\nu}\tilde{u}_{\nu}|^2\Big\}\ds_{\rho}
\leq C\,,
\:\:\:\:\text{for all}\:\nu\,.
$$
Hence, $\{\|\tilde{\eta}_{\nu}\tilde{u}_{\nu}\|_{H^1(\d^+B^+_{\rho}(z_0))}\}$ is uniformly bounded and, since the embedding  
$$
H^1(\d^+B^+_{\rho}(z_0))\subset H^{1/2}(\d^+B^+_{\rho}(z_0))
$$ 
is compact, we can assume that
$$
\tilde{\eta}_{\nu}\tilde{u}_{\nu}\to u
\:\:\text{in}\:\:H^{1/2}(\d^+B^+_{\rho}(z_0))\,,
\:\:\:\text{as}\:\nu\to\infty\,.
$$

We set $\A=B^+_{3r}(z_0)-\overline{B^+_{\rho}(z_0)}$ and let $\{\phi_{\nu}\}\subset D^1(\Rn)$ be such that
$$
\phi_{\nu}=
\begin{cases}
\tilde{\eta}_{\nu}\tilde{u}_{\nu}-u\,,&\text{in}\:B^+_{\rho+\e}(z_0)\,,
\\
0\,,&\text{in}\:\Rn\backslash B^+_{3r-\e}(z_0)\,,
\end{cases}
$$
with $\e>0$ small. Then
$$
\|\tilde{\eta}_{\nu}\tilde{u}_{\nu}-u\|_{H^{1/2}(\d^+B^+_{\rho}(z_0))}
=\|\phi_{\nu}\|_{H^{1/2}(\d^+B^+_{\rho}(z_0))}
\to 0\,,
\:\:\:\:\text{as}\:\nu\to\infty
$$
and there exists $\{\phi_{\nu}^0\}\subset D^1(\A)$ such that
$$
\|\phi_{\nu}+\phi_{\nu}^0\|_{H^1(\A)}
\leq C\|\phi_{\nu}\|_{H^{1/2}(\d^+\A)}
= C\|\phi_{\nu}\|_{H^{1/2}(\d^+B^+_{\rho}(z_0))}
$$
for some $C>0$ independent of $\nu$. Here, $D^1(\A)$ is the closure of $C^{\infty}_0(\A)$ in $H^1(\A)$ and we have set $\d^+\A=\d \A\cap(\Rn\backslash\d\Rn)$ and $\d '\A=\d\A\cap\d\Rn$.

The sequence of functions $\{\zetaup_{\nu}\}=\{\phi_{\nu}+\phi_{\nu}^0\}\subset D^1(\Rn)$ satisfies
$$
\zetaup_{\nu}=
\begin{cases}
\tilde{\eta}_{\nu}\tilde{u}_{\nu}-u\,,&\text{in}\:\overline{B^+_{\rho}(z_0)}\,,
\\
\phi_{\nu}+\phi_{\nu}^0\,,&\text{in}\:B^+_{3r}(z_0)\backslash \overline{B^+_{\rho}(z_0)}\,,
\\
0\,,&\text{in}\:\Rn\backslash B^+_{3r}(z_0)\,.
\end{cases}
$$
In particular, $\zetaup_{\nu}\to 0$ in $H^1(\A)$. 
We set
$$
\tilde{\zetaup}_{\nu}(x)=R_{\nu}^{\frac{n-2}{2}}\zetaup_{\nu}(R_{\nu}\psi_{\nu}^{-1}(x))\,,
\:\:\:\:\text{if}\:\:x\in\psi_{\nu}(B^+_{6r_1}(0))\,,
$$
and $\tilde{\zetaup}_{\nu}(x)=0$ otherwise. Since we are assuming $|z_0|<3(2R-1)r_1\leq 3(2R_{\nu}-1)r_1$ for all $\nu$ large, $\displaystyle{B^+_{3r}(z_0)\subset B_{6r_1R_{\nu}}^+(0)}$. Hence, 
$$
\tilde{\zetaup}_{\nu}(x)=
\begin{cases}
R_{\nu}^{\frac{n-2}{2}}(\tilde{\eta}_{\nu}\tilde{u}_{\nu}-u)(R_{\nu}\psi_{\nu}^{-1}(x))
\,,&\text{if}\:\:x\in\psi_{\nu}(R_{\nu}^{-1}\overline{B^+_{\rho}(z_0)})\,,
\\
R_{\nu}^{\frac{n-2}{2}}(\phi_{\nu}+\phi_{\nu}^0)(R_{\nu}\psi_{\nu}^{-1}(x))\,,&\text{if}\:\:x\in\psi_{\nu}\big{(}R_{\nu}^{-1}(\overline{B^+_{3r}(z_0)}\backslash B^+_{\rho}(z_0))\big{)}\,,
\end{cases}
$$
and $\tilde{\zetaup}_{\nu}(x)=0$ otherwise, and
\ba\label{step4:6}
dI_g(\hat{u}_{\nu})\cdot\tilde{\zetaup}_{\nu}
&=dI_g(\hat{\eta}_{\nu}\hat{u}_{\nu})\cdot\tilde{\zetaup}_{\nu}
\\
&=\int_{B^+_{3r}(z_0)}<d(\tilde{\eta}_{\nu}\tilde{u}_{\nu}),d\zetaup_{\nu}>_{\tg}\dv_{\tg}
-\int_{\d 'B^+_{3r}(z_0)}
|\tilde{\eta}_{\nu}\tilde{u}_{\nu}|^{\frac{2}{n-2}}(\tilde{\eta}_{\nu}\tilde{u}_{\nu})\zetaup_{\nu}
\ds_{\tg}\,,\notag
\end{align}
where $\hat{\eta}_{\nu}(x)=\tilde{\eta}(r_0^{-1}\psi_{\nu}^{-1}(x))$.

Since there exists $C>0$ such that 
$\|\tilde{\zetaup}_{\nu}\|_{H^1(M)}\leq C\|\zetaup_{\nu}\|_{D^1(\Rn)}$, the sequence $\{\tilde{\zetaup}_{\nu}\}$ is uniformly bounded in $H^1(M)$. Hence,
\begin{equation}\label{step4:7}
dI_g(\hat{u}_{\nu})\cdot \tilde{\zetaup}_{\nu}\to 0
\:\:\:\:\text{as}\:\:\nu\to\infty\,.
\end{equation}
Noting that $\zetaup_{\nu}\to 0$ in $H^1(\A)$ and $\zetaup_{\nu}\rightharpoonup 0$ in $D^1(\Rn)$, we obtain
\ba\label{step4:8}
\int_{B^+_{3r}(z_0)}<d(\tilde{\eta}_{\nu}\tilde{u}_{\nu}),d\zetaup_{\nu}>_{\tg}\dv_{\tg}
&= \int_{B^+_{\rho}(z_0)}<d(\zetaup_{\nu}+u),d\zetaup_{\nu}>_{\tg}\dv_{\tg} + o(1)
\\
&=\int_{\Rn}|d\zetaup_{\nu}|_{\tg}^2\dv_{\tg} + o(1)\,.\notag
\end{align}

Similarly,
\begin{equation}\label{step4:9}
\int_{\d 'B_{3r}^+(z_0)}
|\tilde{\eta}_{\nu}\tilde{u}_{\nu}|^{\frac{2}{n-2}}(\tilde{\eta}_{\nu}\tilde{u}_{\nu})\zetaup_{\nu}
\,\ds_{\tg}
=\int_{\d\Rn}|\zetaup_{\nu}|^{\critbordo}\ds_{\tg} + o(1)\,.
\end{equation}
Using (\ref{step4:6}), (\ref{step4:7}), (\ref{step4:8}) and (\ref{step4:9}) we conclude that
\begin{equation}\label{step4:10}
\int_{\Rn}|d\zetaup_{\nu}|_{\tg}^2\dv_{\tg}
=\int_{\d\Rn}|\zetaup_{\nu}|^{\critbordo}\ds_{\tg}+o(1)\,.
\end{equation}
Using again the facts that $\zetaup_{\nu}\to 0$ in $H^1(\A)$ and $\zetaup_{\nu}\rightharpoonup 0$ in $D^1(\Rn)$, we can apply the inequality
$$
\big| |\tilde{\eta}_{\nu}\tilde{u}_{\nu}-u|^{\critbordo}-|\tilde{\eta}_{\nu}\tilde{u}_{\nu}|^{\critbordo}+|u|^{\critbordo}\big|
\leq 
C|u|^{\frac{n}{n-2}}|\tilde{\eta}_{\nu}\tilde{u}_{\nu}-u|
+C|\tilde{\eta}_{\nu}\tilde{u}_{\nu}-u|^{\frac{n}{n-2}}|u|
$$
to see that
$$
\int_{\d\Rn}|\zetaup_{\nu}|^{\critbordo}\ds_{\tg}
=\int_{\d 'B^+_{\rho}(z_0)}|\tilde{\eta}_{\nu}\tilde{u}_{\nu}|^{\critbordo}\ds_{\tg}
-\int_{\d 'B^+_{\rho}(z_0)}|u|^{\critbordo}\ds_{\tg}+o(1)\,.\notag
$$
This implies
\ba\label{step4:10'}
\int_{\d\Rn}|\zetaup_{\nu}|^{\critbordo}\ds_{\tg}
&\leq
\int_{\d 'B^+_{\rho}(z_0)}|\tilde{\eta}_{\nu}\tilde{u}_{\nu}|^{\critbordo}\ds_{\tg}+o(1)
\\
&=
\int_{\d 'B^+_{\rho}(z_0)}|\tilde{u}_{\nu}|^{\critbordo}\ds_{\tg}+o(1)\,,\notag
\end{align}
where we have used the fact that $\tilde{\eta}_{\nu}(z)=1$ for all $z\in B^+_{\rho}(z_0)$.

If $N=N(n)\in\mathbb{N}$ is such that $\partial 'B^+_2(0)$ is covered by $N$ discs in $\partial \mathbb{R}^n_+$ of radius $1$ with center in $\partial 'B^+_2(0)$, then we can choose points $z_i\in\partial ' B^+_{2r}(z_0)$, $i=1,...,N$, satisfying 
$$
\d 'B^+_{\rho}(z_0)\subset\d 'B^+_{2r}(z_0)\subset \bigcup_{i=1}^{N}\d 'B^+_r(z_i)\,.
$$
Hence, using (\ref{step4:3}), (\ref{step4:10}) and (\ref{step4:10'}) we see that
\begin{equation}\label{step4:11}
\int_{\Rn}|d\zetaup_{\nu}|_{\tg}^2\dv_{\tg}+o(1)
= \int_{\d\Rn}|\zetaup_{\nu}|^{\critbordo}\ds_{\tg}
\leq N\l+o(1)\,.
\end{equation}

It follows from (\ref{step4:5}), (\ref{step4:5'}) and the Sobolev inequality (\ref{des:sobolev}) that 
\ba
\left(\int_{\d\Rn}|\zetaup_{\nu}|^{\critbordo}\ds_{\tg}\right)^{\frac{n-2}{n-1}}
&\leq
2^{\frac{n-2}{n-1}}\left(\int_{\d\Rn}|\zetaup_{\nu}|^{\critbordo}dx\right)^{\frac{n-2}{n-1}}\notag
\\
&\leq 2^{\frac{n-2}{n-1}}K_n^2\int_{\Rn}|d\zetaup_{\nu}|^2dx
\leq 2^{1+\frac{n-2}{n-1}}K_n^2\int_{\Rn}|d\zetaup_{\nu}|_{\tg}^2\dv_{\tg}\,.\notag
\end{align}
Then using (\ref{step4:10}) and (\ref{step4:11}) we obtain
\ba
\int_{\Rn}|d\zetaup_{\nu}|_{\tg}^2\dv_{\tg}
&=
\int_{\d\Rn}|\zetaup_{\nu}|^{\critbordo}\ds_{\tg}+o(1)\notag
\\
&\leq
\left( 2^{1+\frac{n-2}{n-1}}K_n^2\right)^{\frac{n-1}{n-2}}
\left(\int_{\Rn}|d\zetaup_{\nu}|_{\tg}^2\dv_{\tg}\right)^{\frac{n-1}{n-2}}+o(1)\notag
\\
&\leq
2^{1+\frac{n-1}{n-2}}K_n^{\frac{2(n-1)}{n-2}}(N\l+o(1))^{\frac{1}{n-2}}
\int_{\Rn}|d\zetaup_{\nu}|_{\tg}^2\dv_{\tg}+o(1)\,.\notag
\end{align}
Now we set $\displaystyle\l_1=\frac{K_n^{-2(n-1)}}{2^{2n-3}N}$ and assume that $\l<\l_1$. Then 
$$
2^{1+\frac{n-1}{n-2}}(N\l)^{\frac{1}{n-2}}K_n^{\frac{2(n-1)}{n-2}}<1,
$$ 
and we conclude that
$$
\lim_{\nu\to\infty}\int_{\Rn}|d\zetaup_{\nu}|_{\tg}^2\dv_{\tg}=0\,.
$$
Hence, $\zetaup_{\nu}\to 0$ in $D^1(\Rn)$. Since $r\leq \rho$, we have
\begin{equation}\label{step4:12}
\tilde{\eta}_{\nu}\tilde{u}_{\nu}\to u
\:\:\:\:\text{in}\:\:H^1(B^+_r(z_0))\,.
\end{equation}

Now let us choose any $z_0=((z_0)^1,...,(z_0)^n)\in \Rn$ satisfying $(z_0)^n>\displaystyle\frac{r}{2}$ and $|z_0|<3(2R-1)r_1$. Using this choice of $z_0$ and $\displaystyle r'=\frac{r}{6}$ replacing $r$, the process above can be performed with some  obvious modifications. In this case, we have $\d 'B^+_{3r'}(z_0)=\emptyset$ and the boundary integrals vanish. Hence, the equality (\ref{step4:10}) already implies that $\tilde{\eta}_{\nu}\tilde{u}_{\nu}\to u$ in $H^1(B^+_{r'}(z_0))$.

If $N_1=N_1(R,n)\in\mathbb{N}$ and $N_2=N_2(R,n)\in\mathbb{N}$ are such that the half-ball $B^+_{2R}(0)$ is covered by $N_1$ half-balls of radius $1$ with center in $\d 'B^+_{2R}(0)$ plus $N_2$ balls of radius $1/6$ with center in $\{z=(z^1,...,z^n)\in 'B^+_{2R}(0)\,;\:z^n>1/2\}$, then the half-ball $B^+_{2Rr}(0)$ is covered by $N_1$ half-balls of radius $r$ with center in $\d 'B^+_{2Rr}(0)$ plus $N_2$ balls of radius $r/6$ with center in $\{z=(z^1,...,z^n)\in B^+_{2Rr}(0)\,;\:z^n>r/2\}$. 

Hence,
$\tilde{\eta}_{\nu}\tilde{u}_{\nu}\to u$ in $H^1(B^+_{2Rr}(0))$, finishing the proof of Claim 1.

\bigskip
Using (\ref{step4:4}),  (\ref{step4:5'}) and Claim 1 with $R=1$ we see that
\ba
\l=
\int_{\d 'B^+_{r}(0)}|\tilde{u}_{\nu}|^{\critbordo}\ds_{\tg}
=\int_{\d 'B^+_{r}(0)}|\tilde{\eta}_{\nu}\tilde{u}_{\nu}|^{\critbordo}\ds_{\tg}
\\
\leq 
2\int_{\d 'B^+_{r}(0)}|u|^{\critbordo}dx+o(1)\,.\notag
\end{align}
It follows that $u\nequiv 0$, due to (\ref{des:sobolev}).

\bigskip\noindent
{\it{Claim 2.}} We have $\lim_{\nu\to\infty}R_{\nu}=\infty$. In particular, Claim 1 can be stated for any $R\geq 1$.

\vspace{0.2cm}\noindent
{\it{Proof of Claim 2.}} Suppose by contradiction that, up to a subsequence, $R_{\nu}\to R'$ as $\nu\to\infty$, for some $1\leq R'<\infty$. Then, since $\u\rightharpoonup 0$ in $H^1(M)$, we have $\tilde{u}_{\nu}\rightharpoonup 0$ in $H^1(B^+_{2r}(0))$. This contradicts the fact that 
$$
\tilde{u}_{\nu}\tilde{\eta}_{\nu}\to u\nequiv 0\,,
\:\:\:\:\text{in}\:\:H^1(B^+_{2r}(0))\,,
$$
which is obtained by applying Claim 1 with $R=1$. This proves Claim 2.

\bigskip
That $u$ is a (weak) solution of (\ref{eq:u}) follows easily from the fact that $\{\hat{u}_{\nu}\}$ is a Palais-Smale sequence for $I_g$ and $\tilde{\eta}_{\nu}\tilde{u}_{\nu}\to u$ in $D^1(\Rn)$.

Now we set 
$$
V_{\nu}(x)=\eta_{\nu}(x)R_{\nu}^{\frac{n-2}{2}}u(R_{\nu}\psi_{x_{\nu}}^{-1}(x))
$$
for $x\in\psi_{x_{\nu}}(B^+_{2r_0}(0))$  and $0$ otherwise.
The proof of the following claim is totally analogous to Step 3 on p.37 of \cite{druet-hebey-robert} with some obvious modifications.

\bigskip\noindent
{\it{Claim 3.}} We have $\hat{u}_{\nu}-V_{\nu}\rightharpoonup 0$, as $\nu\to\infty$, in $H^1(M)$. Moreover, as $\nu\to\infty$,
$$
dI_g(V_{\nu})\to 0 
\:\:\:\text{and}\:\:\:
dI_g(\hat{u}_{\nu}-V_{\nu})\to 0
$$
strongly in $H^1(M)'$, and
$$
I_g(\hat{u}_{\nu})-I_g(\hat{u}_{\nu}-V_{\nu})\to E(u).
$$

\vspace{0.2cm}

We finally observe that if $r_0'>0$ is also sufficiently small then $|(\eta_{\nu}-\eta'_{\nu})V_{\nu}|\to 0$ as $\nu\to \infty$, where $\eta'_{\nu}$ is a smooth cutoff function such that   $\eta'_{\nu}\equiv 1$ in $\psi_{x_{\nu}}(B^+_{r_0'}(0))$ and $\eta'_{\nu}\equiv 0$ in $M\backslash \psi_{x_{\nu}}(B^+_{2r_0'}(0))$.
Hence, the statement of  Lemma \ref{step4} holds for any $r_0>0$ sufficiently small, finishing the proof.
\ep

\bigskip
\bp[Proof of Theorem \ref{struwe_thm}]
According to Lemma \ref{step1}, the Palais-Smale sequence $\{u_{\nu}\}$ for $I_g^{\nu}$ is uniformly bounded in $H^1(M)$. Hence, we can assume that $u_{\nu}\rightharpoonup u^0$ in $H^1(M)$, and $u_{\nu}\to u^0$ a.e in $M$, for some $0\leq u^0\in H^1(M)$. By Lemma \ref{step2}, $u^0$ is a solution to the equations (\ref{eq:Einfty}). Moreover, $\hat{u}_{\nu}=u_{\nu}-u^0$ is Palais-Smale for $I_g$ and satisfies 
$$
I_g(\hat{u}_{\nu})=I_g^{\nu}(u_{\nu})-I_g^{\infty}(u^0)+o(1)\,.
$$

If $\hat{u}_{\nu}\to 0$ in $H^1(M)$, then the theorem is proved. If $\hat{u}_{\nu}\rightharpoonup 0$ in $H^1(M)$ but not strongly, then we apply Lemma \ref{step4} to obtain a new Palais-Smale sequence $\{\hat{u}_{\nu}^1\}$ satisfying 
$$
I_g(\hat{u}_{\nu}^1)\leq I_g(\hat{u}_{\nu})-\beta^*+o(1)=I_g^{\nu}(u_{\nu})-I_g^{\infty}(u^0)-\beta^*+o(1)\,,
$$
where $\beta^*=\displaystyle{\frac{K_n^{-2(n-1)}}{2(n-1)}}$. The term $\beta^*$ appears in the above inequality because $E(u)\geq \beta^*$ for any nontrivial solution $u\in D^1(\Rn)$ to the equations (\ref{eq:Enu}). This can be seen using the Sobolev inequality (\ref{des:sobolev}).

Now we again have either $\hat{u}_{\nu}^1\to 0$ in $H^1(M)$, in which case the theorem is proved, or we apply Lemma \ref{step4} to obtain a new Palais-Smale sequence $\{\hat{u}_{\nu}^2\}$. The process follows by induction and stops by virtue of Lemma \ref{step3}, once we obtain a Palais-Smale sequence $\{\hat{u}_{\nu}^m\}$ with $I_g(\hat{u}_{\nu}^m)$ converging to some $\beta<\beta^*$. 

We are now left with the proof of (\ref{eq:brezis-coron}) and the fact that the $U^j$' obtained by the process above are of the form (\ref{def:U}). To that end, we can follow the proof of Lemma 3.3 in \cite{druet-hebey-robert}, with some simple changes, to obtain the relation (\ref{eq:brezis-coron}) and to prove that  the $U^j$ are nonnegative. For the reader's convenience this is outlined below.

\bigskip\noindent
{\it{Claim.}} The functions $u^0$ and $U^j$ obtained above are nonnegative. Moreover, the identity (\ref{eq:brezis-coron}) holds.

\vspace{0.2cm}\noindent
{\it{Proof of the Claim.}}
That $u^0$ is nonnegative is straightforward. In order to prove that the  $U^j$ are also nonnegative we set $\hat{u}_{\nu}=u_{\nu}-u^0$ and $\displaystyle\mu_{\nu}^j=1/R_{\nu}^j$.

Given integers $N\in [1,m]$ and $s\in [0,N-1]$, we will prove that there exist an integer $p$ and sequences $\{\tilde{x}_{\nu}^k\}_{\nu\in\mathbb{N}}\subset \d M$ and  $\{\lambda_{\nu}^k>0\}_{\nu\in\mathbb{N}}$, for each $k=1,...,p$, such that $d_g(x_{\nu}^N,\tilde{x}_{\nu}^k)/\mu_{\nu}^ N$ is bounded and $\displaystyle \lim_{\nu\to\infty}\lambda_{\nu}^k/\mu_{\nu}^ N=0$, and such that 
\begin{equation}\label{positiv:1}
\int_{\Omega_{\nu}^ N(R)\backslash \bigcup_{k=1}^{p}\tilde{\Omega}_{\nu}^ k(R')}\Big|\hat{u}_{\nu}-\sum_{j=1}^{s}u_{\nu}^ j-u_{\nu}^ N\Big|^{\frac{2n}{n-2}}\dv_g=o(1)+\e(R')\,,
\end{equation}
for any $R,R'>0$. Here, $\Omega_{\nu}^ N(R)=\psi_{x_{\nu}^N}(B^+_{R\mu_{\nu}^N}(0))$,
$\tilde{\Omega}_{\nu}^ k(R')=\psi_{\tilde{x}_{\nu}^k}(B^+_{R'\lambda_{\nu}^k}(0))$ and $\e(R')\to 0$ as $R'\to\infty$. 

We prove (\ref{positiv:1}) by reverse induction on $s$. It follows from Claim 2 in the proof of Lemma \ref{step4} that
$$
\int_{\Omega_{\nu}^ N(R)}\Big|\hat{u}_{\nu}-\sum_{j=1}^{N-1}u_{\nu}^ j-u_{\nu}^ N\Big|^{\frac{2n}{n-2}}\dv_g=o(1)\,,
$$
so that (\ref{positiv:1}) holds for $s=N-1$.

Assuming (\ref{positiv:1}) holds for some $s\in[1,N-1]$, let us prove it does for $s-1$.  

We first consider the case when $d_g(x_{\nu}^s,x_{\nu}^N)$ does not converge to zero as $\nu\to\infty$. In this case, we can assume $\Omega_{\nu}^N(R)\cap\Omega_{\nu}^s(\tilde{R})=\emptyset$ for any $\tilde{R}>0$. Then after rescaling we have
\begin{equation}\label{positiv:2}
\int_{\Omega_{\nu}^ N(R)\backslash \bigcup_{k=1}^{p}\tilde{\Omega}_{\nu}^ k(R')}|u_{\nu}^ s|^{\frac{2n}{n-2}}\dv_g\leq C\int_{\Rn\backslash B^+_{\tilde{R}}(0)}|U^ s|^{\frac{2n}{n-2}}dy.
\end{equation}
Since $\tilde{R}>0$ is arbitrary and $U^ s\in L^{\frac{2n}{n-2}}(\Rn)$, the left side of (\ref{positiv:2}) converges to zero as $\nu\to\infty$. Hence, (\ref{positiv:1}) still holds replacing $s$ by $s-1$. 

Let us now consider the case when $d_g(x_{\nu}^s,x_{\nu}^N)\to 0$ as $\nu\to\infty$. According to Claim 2 in the proof of Lemma \ref{step4}, given $\tilde{R}>0$  we have
$$
\int_{\Omega_{\nu}^ s(\tilde{R})}\Big|\hat{u}_{\nu}-\sum_{j=1}^{s}u_{\nu}^ j\Big|^{\frac{2n}{n-2}}\dv_g=o(1)\,.
$$
Using the induction hypothesis  (\ref{positiv:1}) we then conclude that 
$$
\int_{(\Omega_{\nu}^ N(R)\backslash \bigcup_{k=1}^{p}
\tilde{\Omega}_{\nu}^ k(R'))\cap \Omega_{\nu}^ s(\tilde{R})}
|u_{\nu}^ N|^{\frac{2n}{n-2}}\dv_g=o(1)+\e(R')\,.
$$
First assume that $d_g(x_{\nu}^ s, x_{\nu}^ N)/\mu_{\nu}^N\to\infty$. Rescaling by $\mu_{\nu}^ N$ and using coordinates centered at  $x_{\nu}^ N$, it's not difficult to see that $d_g(x_{\nu}^ s, x_{\nu}^ N)/\mu_{\nu}^s\to\infty$. Hence we can assume that $\Omega_{\nu}^N(R)\cap\Omega_{\nu}^s(\tilde{R})=\emptyset$ for any $\tilde{R}>0$ and we proceed as in the case when  $d_g(x_{\nu}^ s, x_{\nu}^ N)$ does not converge to $0$ to conclude that (\ref{positiv:1}) holds for $s-1$.

If $d_g(x_{\nu}^ s, x_{\nu}^ N)/\mu_{\nu}^N$ does not go to infinity, we can assume that it converges. In this case one can check that $\mu_{\nu}^s/\mu_{\nu}^N\to 0$. We set $\tilde{x}_{\nu}^{p+1}=x_{\nu}^ s$ and $\lambda_{\nu}^{p+1}=\mu_{\nu}^s$, so that $\lambda_{\nu}^{p+1}/\mu_{\nu}^N\to 0$ as $\nu\to\infty$. Observing that  
$$
\int_{\Omega_{\nu}^ N(R)\backslash \bigcup_{k=1}^{p+1}
\tilde{\Omega}_{\nu}^ k(R')}|u_{\nu}^ s|^{\frac{2n}{n-2}}\dv_g
\leq
\int_{M\backslash \Omega_{\nu}^ s(R')}|u_{\nu}^ s|^{\frac{2n}{n-2}}\dv_g
\leq\e(R')\,,
$$
it follows that (\ref{positiv:1}) holds when we replace $p$ by $p+1$ and $s$ by $s-1$.

This proves (\ref{positiv:1}). The above also proves (\ref{eq:brezis-coron}).

We fix an integer $N\in [1,m]$ and $s=0$. Let $\tilde{y}_{\nu}^ k\in\partial\mathbb{R}^{n}_+$ be such that $\tilde{x}_{\nu}^k=\psi_{x_{\nu}^N}^ N(\mu_{\nu}^ N\tilde{y}_{\nu}^ k)$, for $k=1,...,p$. For each $k$, the sequence $\{\tilde{y}^k_{\nu}\}_{\nu\in\mathbb{N}}$ is bounded so there exists $\tilde{y}^k\in\partial\mathbb{R}^{n}_+$ such that $\lim_{\nu\to\infty}\tilde{y}^k_{\nu}=\tilde{y}^k$, possibly after taking a subsequence. Let us set $\tilde{X}=\bigcup_{k=1}^{p}\tilde{y}^k$ and
$$
\tilde{u}_{\nu}^N(y)=(\mu_{\nu}^ N)^{\frac{n-2}{2}}\hat{u}_{\nu}^N(\psi_{x_{\nu}^ N}(\mu_{\nu}^ N y))\,.
$$
It follows from (\ref{positiv:1}) that 
$$
\tilde{u}_{\nu}^N\to U^N\,,
\:\:\:
\text{in}\:L_{loc}^{\frac{2n}{n-2}}(B^+_R(0)\backslash \tilde{X})\,,\:\:\text{as}\:\nu\to\infty\,.
$$
Therefore we can assume that $\tilde{u}_{\nu}\to U^N$ a.e. in $\Rn$ as $\nu\to\infty$. 

If we set
$$
\tilde{u}_{\nu}^{0,N}(y)=(\mu_{\nu}^ N)^{\frac{n-2}{2}}u^ 0(\psi_{x_{\nu}^ N}(\mu_{\nu}^ N y)),
$$
it's easy to prove that 
$$
\tilde{u}_{\nu}^{0,N}\to 0\,,
\:\:\:
\text{in}\:L_{loc}^{\frac{2n}{n-2}}(B^+_R(0))\,,\:\:\text{as}\:\nu\to\infty\,.
$$
Hence, $\tilde{u}_{\nu}^{0,N}\to 0$ a.e. in $\Rn$ as $\nu\to\infty$. Setting 
$$
v_{\nu}^N(y)=(\mu_{\nu}^ N)^{\frac{n-2}{2}}u_{\nu}^N(\psi_{x_{\nu}^ N}(\mu_{\nu}^ N y)),
$$
we see that $v_{\nu}^N\to U^N$ a.e. in $\Rn$ as $\nu\to\infty$. In particular, $U^N$ is nonnegative.
This proves the Claim.

\begin{remark}\label{rm:reg}
For the regularity of the  $U^j$ we can use \cite[Th\'{e}or\`{e}me 1]{cherrier}. Although this theorem is established for compact manifolds we can use the conformal equivalence between $\Rn$ and $B^n\backslash \{\text{point}\}$ and a removable singularities theorem (see Lemma 2.7 on p.1821 of \cite{almaraz3}) to apply it in $B^n$.
\end{remark}

Thus we are able to use the result in \cite{li-zhu} to conclude that the $U^j$ are of the form (\ref{def:U}), so we can write $U^j=U_{\e_j,a_j}$.

This finishes the proof of Theorem \ref{struwe_thm}.
\ep

\bigskip\noindent
{\bf{Acknowledgment.}} I am grateful to the hospitality of the  Department of  Mathematics of the Imperial College London, where part of this work was written during the spring of 2012.


\bigskip\noindent
\textsc{Instituto de Matem\'{a}tica\\Universidade Federal Fluminense\\Niter\'{o}i - RJ, Brazil}
\\{\bf{almaraz@vm.uff.br}}
\end{document}